\documentclass[12pt,a4paper]{article}
\usepackage[scale=.75]{geometry}
\usepackage[affil-it]{authblk}
\usepackage{amssymb,amsmath,amsthm,amsbsy}
\usepackage{physics} 
\usepackage[pdftex,bookmarksnumbered]{hyperref}
\usepackage{graphicx}
\usepackage{epstopdf}
\usepackage{epsfig}
\usepackage{multicol}
\usepackage{subfigure}
\usepackage{multirow}
\usepackage{setspace}
\usepackage{color}
\usepackage[dvipsnames]{xcolor}
\usepackage{lineno}
\usepackage{diagbox}
\usepackage{rotating}

\usepackage{tcolorbox}
\usepackage[width=\textwidth]{caption}
\usepackage{lscape}
\usepackage{booktabs}
\usepackage{adjustbox}

\usepackage[semibold]{ebgaramond}
\usepackage{microtype} 

\newcommand{\trp}{^{\mathsf{T}}}
\usepackage{mathrsfs}

\usepackage{fancyhdr} 
\title{\Large\textbf{A unified perspective of Gaussian process approximation for differential equations}}
\author{Mengwu Guo\thanks{E-mail: \texttt{mengwu.guo@math.lu.se}.}}
\affil{Centre for Mathematical Sciences, Lund University, Sweden}
\date{}

\begin{document}

\maketitle 
\thispagestyle{empty}

\begin{abstract}
\noindent The use of Gaussian processes for approximating differential equations has expanded rapidly, leading to a growing, diverse, and fragmented body of numerical methods. We present a unified Bayesian perspective that places these techniques within a common probabilistic framework, based on a derivative matching interpretation for incorporating differential equation constraints into likelihood. This unified perspective supports both parameter estimation and solution approximation, and shows how a range of existing methods can be understood within it. This work aims to consolidate current developments and provide a foundation for future research.

\end{abstract}

\section{Introduction}

In recent years, there has been a rapidly growing body of work on the use of Gaussian processes \cite{rasmussen2005GPsforMLbook,murphy2012machine,kanagawa2025gaussian} for the approximation of ordinary and partial differential equations, particularly within the emerging field of scientific machine learning. 
These numerical methods seek to combine two complementary sources of information: observational data and physical knowledge expressed through differential equations. On the data-driven side, Gaussian process regression provides a flexible probabilistic representation of the solution function, which allows observations to inform the approximation and associated uncertainty. On the physics-based side, differential equations encode relationships between the solution and its derivatives, thereby providing physical constraints to guide the inference process. Combining these two sources of information allows physical knowledge to improve data-driven approximations, particularly when observations are limited or noisy, while observational data can facilitate the identification of unknown states and parameters in the governing equations.
As a result, a wide range of methodologies has been developed, including various approaches to Gaussian process-based PDE approximation \cite{Simo2011,raissi2017GPsforlinearPDEs,raissi2018gppdes,Yifan2021,Mark2021,pfortner2022physics,beckers2022gaussian,batlle2025error,meng2023sparse,bai2024gaussian}, latent force models \cite{Alvarez2009,Alvarez2013}, and Gaussian process-based approaches to parameter estimation and system identification \cite{yang2021inference,hansen2023learning,ye2024gaussian,McQuarrie2025}. 

Despite this progress, the literature has become increasingly fragmented, with different methods often presented under varied modeling assumptions and algorithmic frameworks. Many of these techniques, however, share common underlying principles: most notably, the preservation of Gaussianity under linear operations \cite{rasmussen2005GPsforMLbook}, and the use of differential equation constraints in Bayesian inference \cite{box1992bayesian}. This raises the question of whether these seemingly disparate methods can be understood within a single, coherent probabilistic framework.

The aim of this note is to provide such a unified perspective. Building on a Bayesian formulation of Gaussian process approximation for differential equations, we show that a broad class of existing methods can be interpreted as instances of a common probabilistic framework. At the same time, the rapid growth of methods on this theme suggests a need for consolidation. While the diversity of approaches has been valuable in investigating different modeling choices and formulations, it is increasingly important to clarify their relationships and underlying assumptions. In this spirit, this note is intended both as a technical contribution and as a pedagogical overview, aiming to synthesize existing ideas and provide a coherent foundation for future developments.

To formalize the discussion, we consider differential equations of the general form:
\begin{equation}\label{DE}
	\mathscr{L}^{\vb* \varphi} u (s)= g^{\vb*\theta}(u (s)) + f(s)\,,\quad s\in D \,.
\end{equation}
Here, $u$ denotes the (scalar-valued) solution function defined on a domain $D$; the operator $\mathscr{L}^{\vb* \varphi}$ is a linear differential operator parameterized by a vector $\vb* \varphi$; the term $g^{\vb*\theta}\circ u$ represents a (nonlinear) source term that depends on the solution state $u$ and is parameterized by a vector $\vb*\theta$; finally, $f$ is a given forcing term defined on $D$. This formulation covers a wide range of differential problems, including linear and nonlinear equations, and accommodates both forward and inverse problems. The objective is to infer the solution $u$, the parameters $(\vb*\varphi, \vb*\theta)$, or both, from both observational data on the solution and the governing equation.

The central idea of this work is to formulate Gaussian process approximation of differential equations within a unified Bayesian framework, rather than to introduce a new numerical algorithm. By placing a Gaussian process prior on the solution function and constructing likelihood terms that encode both observational data and differential equation constraints, with the latter incorporated through a derivative matching mechanism, we obtain a posterior distribution that jointly describes the solution and the unknown parameters. Parameter estimation and solution approximation can then be addressed as two inferential scenarios through appropriate marginalization of this joint posterior. The same perspective is further extended to weak formulations of differential equations.

The remainder of the paper is organized as follows. In Section 2, we develop the unified Bayesian framework, which gives rise to two typical scenarios discussed in Section 3: parameter estimation and solution approximation. In Section 4, we illustrate the generality of the framework by relating it to several existing methods and extending it from strong-form collocation to weak-form functional evaluation. Finally, conclusions are drawn in Section 5.

For clarity of notation, italic bold symbols are used throughout the paper to denote vector-valued (random) variables, such as the parameters $\vb*\varphi$ and $\vb*\theta$, the latent state $\vb*U$, and the observables $\vb*Y$ and $\vb*R$. Upright bold symbols are adopted to denote vectors and matrices of given data, such as the data vectors $\vb y$ and $\vb f$, and the covariance matrix $\vb K$. We write $p_{\vb*X}(\vb* x)$ for the probability density function (PDF) of a random variable $\vb*X$, where $\vb*x$ is the argument of the PDF. Conditional density functions are written as $p_{\vb*X_1\mid \vb*X_2}(\vb*x_1 \mid \vb* x_2)$, which denotes the PDF of $\vb*X_1$ evaluated at $\vb*x_1$, under the condition $\vb*X_2 = \vb*x_2$.  Following standard convention, we use capitalized and lowercase letters to distinguish random latent variables or observables from their corresponding PDF variables or observed data. For example, $\vb*U$ denotes a latent random vector, whereas $\vb*u$ the corresponding argument in a PDF; similarly, $\vb*Y$ denotes an observable, whereas $\vb y$ its observed value. This distinction is not made for the parameters $\vb*\varphi$ and $\vb*\theta$.

\section{Unified Bayesian framework}

In this section, we present a unified Bayesian formulation for Gaussian process approximation of differential equations, including the specification of prior, likelihood, and posterior distributions. 

\subsection{Prior}

Assume that the solution function $u$ follows a prior distribution given by a Gaussian process on the domain $D$:
\begin{equation}
	u(\cdot) \sim \mathcal{GP}(0, \kappa(\cdot,\cdot))\,,
\end{equation}
which is assumed to have zero mean and is assigned a covariance (kernel) function $\kappa$.  Let $\vb*U$ denote the random vector obtained by evaluating the solution $u$ at $m$ observed locations/coordinates $\mathcal{S} = \{s_1, \cdots. s_m\}\subset D$, i.e., $\vb*U = (u(s_1), \cdots, u(s_m))\trp$. The density function of the prior distribution of $\vb*U$ is then given by \footnote{Throughout this paper, $\mathcal{N}(\cdot \mid \vb a, \vb A)$ denotes the Gaussian density function with mean vector $\vb a$ and covariance matrix $\vb A$, while $\mathcal{N}(\vb a, \vb A)$ denotes the Gaussian distribution.}
\begin{equation}
	p_{\vb*U}(\vb*u) = \mathcal{N}(\vb*{u}\mid \vb 0, \vb K)\,,
\end{equation}
where the covariance matrix $\vb K = \kappa(\mathcal{S},\mathcal{S})$. i.e., $\vb K_{ij} = \kappa(s_i , s_j)$.

Overall, the joint prior distribution of the parameters $\vb* \varphi$ and $\vb* \theta$, together with the state $\vb* U$, is given by
\begin{equation}
	\pi_\text{prior} (\vb* \varphi, \vb* \theta, \vb* u) = p_{\vb* \varphi}(\vb* \varphi) ~ p_{\vb* \theta}(\vb* \theta) ~ p_{\vb*U}(\vb*u)\,.
\end{equation}
Here, $p_{\vb* \varphi}(\vb* \varphi)$ and $p_{\vb* \theta}(\vb* \theta)$ denote prior distributions on the parameters $\vb* \varphi$ and $\vb* \theta$, respectively. These priors are left unspecified here and will be defined in specific cases. We assume that $\vb* \varphi$, $\vb* \theta$, and $\vb* U$ are mutually independent under the prior.

\subsection{Likelihood}

\subsubsection{Observables and likelihood function}

To incorporate both observational data and the differential equation constraint \eqref{DE}, we consider two observables in the definition of the likelihood. The first observable $\vb* Y$ represents a direct observation on the state $\vb*{U}$, subject to additive white noise with variance $\sigma_y^2$, i.e.,
\begin{equation}\label{like_y}
	p_{\vb*Y \mid \vb*U} (\vb y \mid \vb*u) = \mathcal{N}(\vb y\mid \vb*u, \sigma_y^2 \vb I_{m})\,,
\end{equation}
where $\vb y$ collects the observed data of $\vb* Y$. The second observable $\vb* R$ corresponds to the evaluation of $\rho(u; \vb* \varphi, \vb*\theta) := \mathscr{L}^{\vb* \varphi} u - g^{\vb*\theta}\circ u$ at $m'$ locations $\mathcal{S}' = \{s'_1, \cdots, s'_{m'}\}\subset D$. Although $\rho$ is a nonlinear operator acting on the unknown solution $u$, with parameters $\vb*\varphi$ and $\vb*\theta$ yet to be inferred, its output must equal $f$ to satisfy the differential equation \eqref{DE}. Therefore, the given forcing term $f$ provides data $\vb f = (f(s'_1), \cdots, f(s'_{m'}))\trp$, which can be viewed as observations of $\vb* R$. From a mechanics viewpoint, $\rho$ can be seen as the `internal force' of the system, which achieves equilibrium with the `external force' $f$, i.e., $\rho = f$.

A likelihood function that accounts for both observables is given by
\begin{subequations}
\begin{align}
	\pi_\text{like} (\vb y, \vb f \mid \vb* \varphi, \vb* \theta, \vb* u)  & = p_{\vb*Y, \vb*R ~\mid~ \vb* \varphi, \vb* \theta, \vb*U} (\vb y, \vb f \mid \vb* \varphi, \vb* \theta, \vb* u) \\
	& = p_{\vb*Y \mid \vb*U} (\vb y \mid  \vb* u) ~ p_{\vb*R \mid \vb* \varphi, \vb* \theta, \vb*U} (\vb f \mid \vb* \varphi, \vb* \theta, \vb* u)
\end{align}
\end{subequations}
where the first term $p_{\vb*Y \mid \vb*U} (\vb y \mid  \vb* u)$ is given in \eqref{like_y}. The factorization follows from the independence between $\vb*Y \mid \vb*U$ and $\vb*R \mid \vb* \varphi, \vb* \theta, \vb*U$, induced by the additive white noise in the observation model for $\vb*Y$.

\subsubsection{Derivative matching}

Let the random vector $\vb*Z$ denote the evaluation of $\mathscr{L}^{\vb* \varphi} u$ at the locations in $\mathcal{S}'$. We refer to $\vb*Z$ as `derivatives', as it is obtained from the solution function $u$ through differential operations. Then, the likelihood term $p_{\vb*R \mid \vb* \varphi, \vb* \theta, \vb*U} (\vb f \mid \vb* \varphi, \vb* \theta, \vb* u)$ admits the following factorization:
\begin{equation}\label{deriv_match}
		p_{\vb*R \mid \vb* \varphi, \vb* \theta, \vb*U} (\vb f \mid \vb* \varphi, \vb* \theta, \vb* u) = \int p_{\vb*R\mid \vb*Z, \vb*\theta} (\vb f\mid \vb*z, \vb*\theta)~ p_{\vb*Z\mid \vb*\varphi, \vb*U} (\vb*z\mid \vb*\varphi, \vb*u) ~\dd\vb*z \,.
\end{equation}

\paragraph{Derivatives estimated through Gaussian process} Specifically, $\vb*Z\mid \vb*\varphi, \vb*U$ represents a Gaussian process-based estimate of $\vb*Z$ obtained by applying the linear differential operator $\mathscr{L}^{\vb* \varphi}$, given by the conditional Gaussian formula 
\footnote{Consider a joint Gaussian distribution $p_{\vb*X_1, \vb*X_2} (\vb*x_1, \vb*x_2) = \mathcal N \left(
\begin{pmatrix}
	\vb*x_1 \\ \vb*x_2
\end{pmatrix}	\Big |
\begin{pmatrix}
	\vb*\mu_1 \\ \vb*\mu_2
\end{pmatrix}\,,
\begin{bmatrix}
	\vb*\Sigma_{11} & \vb*\Sigma_{12} \\
	\vb*\Sigma_{21} & \vb*\Sigma_{22}
\end{bmatrix}
	\right)$ and assume that $\vb*\Sigma_{22}$ is invertible, then $p_{\vb*X_1 \mid \vb*X_2} (\vb*x_1 \mid \vb*x_2) = \mathcal N\left(\vb*x_1 \mid \vb*\mu_1 + \vb*\Sigma_{12}\vb*\Sigma_{22}^{-1}(\vb*x_2 - \vb*\mu_2), \vb*\Sigma_{11} - \vb*\Sigma_{12}\vb*\Sigma_{22}^{-1} \vb*\Sigma_{21}\right)$.}
as follows:
\begin{equation}
	\begin{split}
		p_{\vb*Z\mid \vb*\varphi, \vb*U} & (\vb*z\mid \vb*\varphi, \vb*u) = \mathcal{N}(\vb*z\mid  (\mathscr{L}^{\vb* \varphi}\otimes \text{Id})\kappa(\mathcal{S}',\mathcal{S})\vb{K}^{-1}\vb*u\,,\\
		& (\mathscr{L}^{\vb* \varphi}\otimes \mathscr{L}^{\vb* \varphi})\kappa(\mathcal{S}',\mathcal{S}')- [(\mathscr{L}^{\vb* \varphi}\otimes \text{Id})\kappa(\mathcal{S}',\mathcal{S})]\vb{K}^{-1} [(\text{Id} \otimes \mathscr{L}^{\vb* \varphi} )\kappa(\mathcal{S},\mathcal{S}')])\,,
	\end{split}
\end{equation}
where $\text{Id}$ denotes the identity operator. Here, the tensor product notation $\mathscr{L}^{\vb* \varphi}\otimes \text{Id}$ and $\mathscr{L}^{\vb* \varphi}\otimes \mathscr{L}^{\vb* \varphi}$ denotes the action of differential operators on the kernel function with respect to its first and second arguments, respectively. For example, $(\mathscr{L}^{\vb* \varphi}\otimes \text{Id})\kappa(s,s')$ means that $\mathscr{L}^{\vb* \varphi}$ is applied to $\kappa(\cdot,s')$ with respect to the first variable, while the second variable is left unchanged. Similarly, $(\mathscr{L}^{\vb* \varphi}\otimes \mathscr{L}^{\vb* \varphi})\kappa(s,s')$ indicates that the operator is applied to both arguments of the kernel. Moreover, expressions such as $\kappa(\mathcal{S}',\mathcal{S})$ and $\kappa(\mathcal{S}',\mathcal{S}')$ denote matrices obtained by evaluating the kernel at pairs of locations in the corresponding sets; for instance, $[\kappa(\mathcal{S}',\mathcal{S})]_{ij} = \kappa(s'_i, s_j)$. The operator-transformed kernel functions reflect the well-known property that \textit{Gaussianity is preserved under linear operations}.

It is worth pointing out that the covariance kernel $\kappa$ is assumed to be sufficiently smooth for the differential operator $\mathscr{L}^{\vb*\varphi}$ acting on its arguments to be well defined. The required regularity depends on the order of the differential operator. Beyond this minimal requirement, the kernel may be selected to reflect the expected smoothness of the solution. For instance, squared-exponential kernels provide high regularity, while Matérn kernels offer adjustable smoothness through their smoothness parameter.

\paragraph{Derivatives from the differential equation constraint} Now consider an estimate of the solution $u$ at the $\mathcal{S}'$ locations (or a direct observation, if available), for example via a Gaussian process regression $\tilde{\vb u} = \kappa(\mathcal{S}',\mathcal{S})(\vb K+\sigma_y^2\vb I)^{-1}\vb y$. The conditional distribution $\vb*R\mid \vb*Z, \vb*\theta$ can then be interpreted as enforcing consistency between the left- and right-hand sides of \eqref{DE} at the locations in $\mathcal{S}'$, given by
\begin{equation}
	p_{\vb*R\mid \vb*Z, \vb*\theta} (\vb f\mid \vb*z, \vb*\theta) = \mathcal{N}(\vb f \mid \vb*z - g^{\vb*\theta}(\tilde{\vb u}), \sigma_f^2\vb I_{m'})\,,
\end{equation}
where additive white noise with variance $\sigma_f^2$ is introduced to model potential misalignment between $\vb*Z - g^{\vb*\theta}(\tilde{\vb u})$ (i.e., an approximation of $\vb*R$) and the observed data $\vb*R = \vb f$ given by the forcing term. This approximation strategy has also been employed in previous Gaussian process-based approaches \cite{yang2021inference,ye2024gaussian,McQuarrie2025}. Here $g^{\vb*\theta}(\tilde{\vb u})$ is a shorthand that stands for the component-wise evaluation of $g^{\vb*\theta}$ on the vector $\tilde{\vb u}$. 

This treatment provides a computationally convenient approximation of the nonlinear source term while retaining a Gaussian structure in the likelihood. Its reliability depends on the accuracy of the solution estimate $\tilde{u}$ and the adequacy of the assumed Gaussian noise model for the discrepancy. If either is inadequate, the resulting inference may be inaccurate. Alternative strategies for handling nonlinearity include lifting techniques, Newton-type methods based on linearization, and additional closure modeling, where appropriate.

\medskip
With both terms in the integrand of \eqref{deriv_match} specified, the integral can be interpreted as a form of derivative matching. The term $p_{\vb*R\mid \vb*Z, \vb*\theta} (\vb f\mid \vb*z, \vb*\theta)$ represents the expression of the derivatives $\vb*Z$ through the right-hand side of the differential equation \eqref{DE}, while $p_{\vb*Z\mid \vb*\varphi, \vb*U} (\vb*z\mid \vb*\varphi, \vb*u)$ provides a Gaussian process-based estimate of these derivatives. These two representations are required to agree, up to Gaussian noise with variance $\sigma_f^2$. The latent variables $\vb*Z$ can therefore be marginalized out, yielding the following explicit expression for the likelihood term $p_{\vb*R \mid \vb* \varphi, \vb* \theta, \vb*U} (\vb f \mid \vb* \varphi, \vb* \theta, \vb* u)$:
\begin{equation}
	\begin{split}
		p_{\vb*R \mid \vb* \varphi, \vb* \theta, \vb*U}  (\vb f & \mid \vb* \varphi, \vb* \theta, \vb* u) = \mathcal{N}(\vb f\mid  (\mathscr{L}^{\vb* \varphi}\otimes \text{Id})\kappa(\mathcal{S}',\mathcal{S})\vb{K}^{-1}\vb*u  - g^{\vb*\theta}(\tilde{\vb u})   \,,\\
		 \sigma_f^2\vb I_{m'}  & +  (\mathscr{L}^{\vb* \varphi}\otimes \mathscr{L}^{\vb* \varphi})\kappa(\mathcal{S}',\mathcal{S}')  - [(\mathscr{L}^{\vb* \varphi}\otimes \text{Id})\kappa(\mathcal{S}',\mathcal{S})]  \vb{K}^{-1} [(\text{Id} \otimes \mathscr{L}^{\vb* \varphi} )\kappa(\mathcal{S},\mathcal{S}')])\,.
	\end{split}
\end{equation}
The derivative-matching formulation is also closely related to the product-of-experts principle \cite{hinton2002training}, which combines multiple probabilistic descriptions of the same quantity. Related Gaussian process-based gradient-matching techniques have been developed for parameter inference in ordinary differential equations \cite{calderhead2008odegps,dondelinger2013ode,tronarp2019probabilisticodes}.

\medskip
This concludes the construction of the likelihood function, which integrates both observational data and model-based constraints.

\subsection{Posterior}

Bayes' rule yields
\begin{equation}\label{post}
	\begin{split}
		\pi_\text{post} (\vb* \varphi, \vb* \theta, \vb* u \mid \vb{y}, \vb{f}) & = p_{\vb* \varphi, \vb* \theta, \vb* U \mid \vb*Y, \vb*R} (\vb* \varphi, \vb* \theta, \vb* u \mid \vb y, \vb f) \\
		& \propto  \pi_\text{like} (\vb y, \vb f \mid \vb* \varphi, \vb* \theta, \vb* u)  ~ \pi_\text{prior} (\vb* \varphi, \vb* \theta, \vb* u) \\
		& = p_{\vb*Y \mid \vb*U} (\vb y \mid  \vb* u) ~ p_{\vb*R \mid \vb* \varphi, \vb* \theta, \vb*U} (\vb f \mid \vb* \varphi, \vb* \theta, \vb* u) \cdot p_{\vb*U}(\vb*u) ~ p_{\vb* \varphi}(\vb* \varphi) ~ p_{\vb* \theta}(\vb* \theta)\,.
	\end{split}
\end{equation}
Note that this posterior distribution involves the reconstructions of the solution values $\vb* U$, not only through denoising the directly observed data $\vb*Y = \vb y$, but also by incorporating the differential equation constraint. In practice, however, this posterior is generally intractable to evaluate in closed form due to the nonlinearity in the model with respect to the parameters $(\vb* \varphi, \vb* \theta)$, and thus requires suitable approximation through appropriate inference techniques.

When using Gaussian processes for approximating differential equations, two typical scenarios may be considered. In the first scenario, the focus is on parameter estimation, i.e., the inference of $(\vb* \varphi, \vb* \theta)$. Conceptually, the posterior distribution of $(\vb* \varphi, \vb* \theta)$ can be obtained by marginalizing out $\vb*U$ in \eqref{post}, i.e., 
\begin{equation}
	\pi_\text{post} (\vb* \varphi, \vb* \theta \mid \vb{y}, \vb{f}) = \int \pi_\text{post} (\vb* \varphi, \vb* \theta, \vb* u \mid \vb{y}, \vb{f})~\dd\vb*u\,.
\end{equation}
However, sampling $\vb*U$ from the posterior is often computationally expensive. Fortunately, this marginalization need not be carried out explicitly at the posterior stage, as both the prior and likelihood can be simplified in this setting. This will be discussed in detail in Section 3.1.

In the second scenario, the goal is to predict the solution $u$ at new, unseen (test) locations using the Gaussian process. To this end, let the random variable $U^* = u(s^*)$ denote the solution value at a test location $s^* \in D$. Its predictive distribution is then given by
\begin{equation}\label{pred}
	p_{U^*\mid \vb*Y, \vb*R} (u^* \mid \vb y, \vb f) = \int p_{U^* \mid \vb*\varphi, \vb*\theta, \vb*U, \vb*Y, \vb*R}(u^*\mid \vb*\varphi, \vb*\theta, \vb*u, \vb y, \vb f) ~ \pi_\text{post} (\vb* \varphi, \vb* \theta, \vb* u \mid \vb{y}, \vb{f})  ~\dd\vb*\varphi ~\dd\vb*\theta ~\dd\vb*u.
\end{equation}
Note that, in fact, $p_{U^* \mid \vb*\varphi, \vb*\theta, \vb*U, \vb*Y, \vb*R}(u^*\mid \vb*\varphi, \vb*\theta, \vb*u, \vb y, \vb f) = p_{U^* \mid \vb*\varphi, \vb*\theta, \vb*U, \vb*R}(u^*\mid \vb*\varphi, \vb*\theta, \vb*u, \vb f)$, i.e., one can remove the condition on $\vb*Y = \vb y$, because $\vb* Y \mid \vb* U$ is defined by an independent white noise. Again, here involves the sampling and marginalization over $\vb*U$, which should be avoided in practice. A computationally efficient strategy will be discussed in Section 3.2.

\begin{figure}
	\centering
	\includegraphics[width=.7\textwidth]{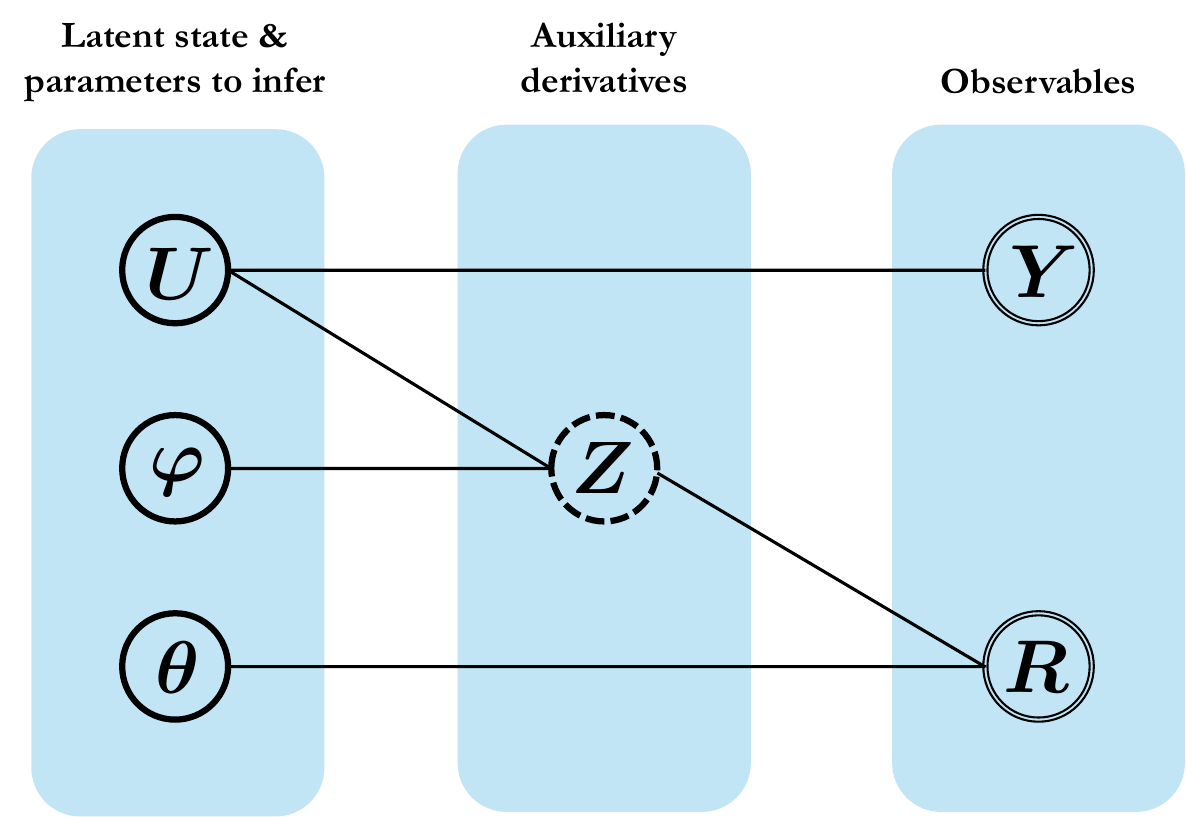}
	\caption{Inferential structure of the unified Bayesian framework}
	\label{fig:ProbStruc}
\end{figure}

\medskip
This completes the Bayesian formulation, bringing together the prior, likelihood, and posterior in a unified framework that forms the basis for the discussions of two typical scenarios in the next section. Furthermore, the overall probabilistic structure of the inference framework is visualized in Figure \ref{fig:ProbStruc}.

\section{Two typical scenarios}

Based on the unified Bayesian formulation developed in Section 2, we consider two typical scenarios that arise in Gaussian process approximation of differential equations, corresponding to parameter estimation and solution approximation.

\subsection{Scenario I: parameter estimation}

When the goal is to estimate parameters $(\vb* \varphi, \vb* \theta)$, the state $\vb*U$ should be marginalized out. This leads to the following integral and factorization, which removes the dependence on $\vb*u$ from the likelihood function:
\begin{equation}
	\pi_\text{like} (\vb{y}, \vb{f} \mid \vb* \varphi, \vb* \theta) =  \int \pi_\text{like} (\vb{y}, \vb{f} \mid \vb* \varphi, \vb* \theta, \vb* u)~p_{\vb*U}(\vb*u) ~\dd\vb*u = p_{\vb*Y } (\vb y ) ~ p_{\vb*R \mid \vb* \varphi, \vb* \theta, \vb*Y} (\vb f \mid \vb* \varphi, \vb* \theta, \vb y)\,.
\end{equation}
One observes that $p_{\vb*Y } (\vb y )  = \mathcal{N}(\vb y \mid \vb 0, \vb K + \sigma_y^2 \vb I_m)$ is independent of $(\vb* \varphi, \vb* \theta) $, and that all dependence on $(\vb* \varphi, \vb* \theta) $ arises from 
\begin{equation}\label{like_f}
	\begin{split}
		& p_{\vb*R \mid \vb* \varphi, \vb* \theta, \vb*Y} (\vb f \mid \vb* \varphi, \vb* \theta, \vb y)  = \mathcal{N}(\vb f\mid  (\mathscr{L}^{\vb* \varphi}\otimes \text{Id})\kappa(\mathcal{S}',\mathcal{S})(\vb{K}+\sigma_y^2 \vb I_m )^{-1}\vb y - g^{\vb*\theta}(\tilde{\vb u})  \,,\\
		& \sigma_f^2\vb I_{m'}  +  (\mathscr{L}^{\vb* \varphi}\otimes \mathscr{L}^{\vb* \varphi})\kappa(\mathcal{S}',\mathcal{S}')   - [(\mathscr{L}^{\vb* \varphi}\otimes \text{Id})\kappa(\mathcal{S}',\mathcal{S})](\vb{K}+\sigma_y^2 \vb I_m )^{-1}  [(\text{Id} \otimes \mathscr{L}^{\vb* \varphi} )\kappa(\mathcal{S},\mathcal{S}')])\,.
	\end{split}
\end{equation}
This expression is obtained by first applying Gaussian process regression to the observations $\vb*Y = \vb y$ in order to estimate the derivatives, and then aligning these estimated derivatives with the right-hand side of the differential equation, in the sense of the derivative matching formulation introduced earlier. Accordingly, the posterior simplifies to
\begin{equation}
	\pi_\text{post} (\vb* \varphi, \vb* \theta \mid \vb{y}, \vb{f}) = p_{\vb* \varphi, \vb* \theta \mid \vb*Y, \vb*R} (\vb* \varphi, \vb* \theta \mid \vb y, \vb f) \propto p_{\vb*R \mid \vb* \varphi, \vb* \theta, \vb*Y} (\vb f \mid \vb* \varphi, \vb* \theta, \vb y) ~ p_{\vb* \varphi}(\vb* \varphi) ~ p_{\vb* \theta}(\vb* \theta)\,.
\end{equation}

To propagate posterior uncertainty in the parameters to the solution, one may view \eqref{DE} as a parametric problem and introduce the parameter-to-solution map
\begin{equation}
	\Psi : (\vb* \varphi, \vb* \theta) \mapsto u(\cdot; \vb* \varphi, \vb* \theta)\,.
\end{equation}
Assuming that $\Psi$ is measurable, the conditional distribution of $u(\cdot)\mid \vb*Y = \vb y, \vb*R = \vb f$ is given by the pushforward measure
\begin{equation}
	\nu_{u(\cdot)\mid \vb*Y = \vb y, \vb*R = \vb f} = \Psi_ \# \nu_{\vb* \varphi, \vb* \theta \mid \vb*Y = \vb y, \vb*R = \vb f}\,,
\end{equation}
where $\nu$ denotes a probability measure, and $\nu_{\vb* \varphi, \vb* \theta \mid \vb*Y = \vb y, \vb*R = \vb f}$ is defined by the density function $\pi_\text{post} (\vb* \varphi, \vb* \theta \mid \vb{y}, \vb{f})$.

\paragraph{Special case: parameter estimation in nonlinear dynamics} (see, e.g., \cite{ye2024gaussian}) Consider nonlinear dynamics described by the ordinary differential equation
\begin{equation}
	\frac{\text{d}}{\text{d}t} u(t) = g^{\vb*\theta}(u(t)) \,,\quad t \in D = (0, T)
\end{equation}
parametrized solely by $\vb*\theta$ through the right-hand side. In this case, the differential operator $\mathscr{L} = \text{d}/\text{d}t$ represents the time derivative and does not involve additional parameters (i.e., no $\vb*\varphi$); the coordinate $s= t$ represents time; the finite sets $\mathcal{S}, \mathcal{S}' \subset D$ collect $m$ and $m'$ time instances, respectively; and $f =0$, i.e., there is no external forcing (or the forcing term is absorbed into $g$). The posterior distribution of the parameter $\vb*\theta$ is then given by
\begin{align}
	\pi_\text{post} (\vb* \theta \mid \vb{y}, \vb{0})  = p_{\vb* \theta \mid \vb*Y, \vb*R} (\vb* \theta \mid \vb y, \vb 0) \propto p_{\vb*R \mid \vb* \theta, \vb*Y} (\vb 0 \mid  \vb* \theta, \vb y)  ~ p_{\vb* \theta}(\vb* \theta) \,.
\end{align}
To simplify the notation in \eqref{like_f}, we introduce
\begin{subequations}
	\begin{align}
		\widehat{\vb*\mu}_z  : =~ & (\mathscr{L}\otimes \text{Id})\kappa(\mathcal{S}',\mathcal{S})(\vb{K}+\sigma_y^2 \vb I_m )^{-1}\vb y\,, \quad \text{and} \\
		\begin{split}
		\widehat{\vb*\Sigma}_z  : =~ & \sigma_f^2\vb I_{m'}  +  (\mathscr{L}\otimes \mathscr{L})\kappa(\mathcal{S}',\mathcal{S}') \\
		&   - [(\mathscr{L}\otimes \text{Id})\kappa(\mathcal{S}',\mathcal{S})](\vb{K}+\sigma_y^2 \vb I_m )^{-1} [(\text{Id} \otimes \mathscr{L})\kappa(\mathcal{S},\mathcal{S}')].
		\end{split}
	\end{align}
\end{subequations}
Here, $\widehat{\vb*\mu}_z$ and $\widehat{\vb*\Sigma}_z$ denote the mean vector and covariance matrix, respectively, for a Gaussian process approximation of the time derivatives $\vb*Z$, conditioned on the state observations $\vb*Y = \vb y$.
The posterior $\pi_\text{post} (\vb* \theta)$ can then be written explicitly as
\begin{equation}
	-\log \pi_\text{post} (\vb* \theta \mid \vb{y}, \vb{0}) = 
	\frac{1}{2}\|g^{\vb*\theta}(\tilde{\vb u}) - \widehat{\vb*\mu}_z \|^2_{\widehat{\vb*\Sigma}_z^{-1}} - \log p_{\vb* \theta}(\vb* \theta) + C\,,
\end{equation}
in which $\|\vb a\|_{\widehat{\vb*\Sigma}_z^{-1}} = (\vb a \trp \widehat{\vb*\Sigma}_z^{-1} \vb a)^{1/2}$ denotes the $\widehat{\vb*\Sigma}_z^{-1}$-weighted norm for $\vb a \in \mathbb{R}^{m'}$, and $C$ is a constant independent of $\vb*{\theta}$. It follows that the maximum a posteriori (MAP) estimate of $\vb*\theta$ is equivalent to a nonlinear, regularized generalized least squares problem. Furthermore, the parameter-to-solution map $\Psi: \vb*\theta \mapsto u(\cdot; \vb*\theta)$ can be approximated by numerical time integration.

An analytical illustration of this parameter estimation scenario is provided in Appendix A, where the joint posterior distribution of the governing parameters in the Lotka-Volterra system is derived in closed form.

\subsection{Scenario II: equation solving via Gaussian process regression}

In this second scenario, the Gaussian process model is used directly within the inference procedure to approximate the solution field. By marginalizing out the state $\vb*U$ first (as discussed in subsection 3.1) and then the parameters, the predictive distribution \eqref{pred} for $U^*= u(s^*) $ at an unseen location $s^*$ can be rewritten as
\begin{equation}
	p_{U^*\mid \vb*Y, \vb*R} (u^* \mid \vb y, \vb f) = \int p_{U^* \mid \vb*\varphi, \vb*\theta, \vb*Y, \vb*R}(u^*\mid \vb*\varphi, \vb*\theta, \vb y, \vb f) ~ \pi_\text{post} (\vb* \varphi, \vb* \theta \mid \vb y, \vb f )  ~\dd\vb*\varphi ~\dd\vb*\theta.
\end{equation}
Here, $U^* \mid \vb*\varphi, \vb*\theta, \vb*Y, \vb*R $ can be interpreted as the prediction of $U^*$ via Gaussian process regression, conditioned on two types of data,  $\vb*Y= \vb y$ and $\vb*R = \vb f$. Since this regression depends on the parameters $(\vb* \varphi, \vb* \theta)$, marginalization requires sampling over their posterior distribution.

To derive the conditional density function
\begin{equation}
	p_{U^* \mid \vb*\varphi, \vb*\theta, \vb*Y, \vb*R}(u^*\mid \vb*\varphi, \vb*\theta, \vb y, \vb f) = \frac{p_{U^* , \vb*Y, \vb*R \mid \vb*\varphi, \vb*\theta }(u^*, \vb y, \vb f\mid \vb*\varphi, \vb*\theta)}{p_{\vb*Y, \vb*R \mid \vb*\varphi, \vb*\theta }(\vb y, \vb f\mid \vb*\varphi, \vb*\theta)}\,,
\end{equation}
one can first write the joint Gaussian density function $p_{U^* , \vb*Y, \vb*R \mid \vb*\varphi, \vb*\theta }(u^*, \vb y, \vb f\mid \vb*\varphi, \vb*\theta)$ and then apply the formula for conditional Gaussian distributions. In particular,
\begin{equation}
	\begin{split}
	&p_{U^* , \vb*Y, \vb*R \mid \vb*\varphi, \vb*\theta, }(u^*, \vb y, \vb f\mid \vb*\varphi, \vb*\theta) \\
	=~ &  \mathcal{N} \left(
	\begin{pmatrix}
		u^* \\ \vb y \\ \vb f
	\end{pmatrix}\, \Bigg |
	\begin{pmatrix}
		0 \\ 0 \\ -g^{\vb*\theta}(\tilde{\vb u})
	\end{pmatrix}\,,  \right.\\
	& ~~~\left.\begin{bmatrix}
		\kappa(s^*, s^*) & \kappa(s^*, \mathcal{S}) & (\text{Id} \otimes \mathscr{L}^{\vb* \varphi} ) \kappa(s^*, \mathcal{S}')\\
		\kappa(\mathcal{S},s^*) & \kappa(\mathcal{S},\mathcal{S})+\sigma_y^2 \vb I_m & (\text{Id} \otimes \mathscr{L}^{\vb* \varphi} ) \kappa(\mathcal{S}, \mathcal{S}')\\
		(\mathscr{L}^{\vb* \varphi}\otimes \text{Id})\kappa(\mathcal{S}',s^*) & (\mathscr{L}^{\vb* \varphi}\otimes \text{Id})\kappa(\mathcal{S}',\mathcal{S}) & (\mathscr{L}^{\vb* \varphi}\otimes \mathscr{L}^{\vb* \varphi})\kappa(\mathcal{S}',\mathcal{S}') + \sigma_f^2 \vb I_{m'}\\
	\end{bmatrix}
	\right)\,.
	\end{split}
\end{equation}
For notational convenience, define
\begin{subequations}
	\begin{align}
		\widehat{\vb* k}^{\vb* \varphi}(s^*) & : = 
		\begin{pmatrix}
			\kappa(\mathcal{S},s^*) \\ (\mathscr{L}^{\vb* \varphi}\otimes \text{Id})\kappa(\mathcal{S}',s^*)
		\end{pmatrix}\,, \quad \text{and} \\
	\widehat{\vb K}^{\vb* \varphi} & : = 
	\begin{bmatrix}
		\kappa(\mathcal{S},\mathcal{S})+\sigma_y^2 \vb I_m & (\text{Id} \otimes \mathscr{L}^{\vb* \varphi} ) \kappa(\mathcal{S}, \mathcal{S}')\\
		(\mathscr{L}^{\vb* \varphi}\otimes \text{Id})\kappa(\mathcal{S}',\mathcal{S}) & (\mathscr{L}^{\vb* \varphi}\otimes \mathscr{L}^{\vb* \varphi})\kappa(\mathcal{S}',\mathcal{S}') + \sigma_f^2 \vb I_{m'}
	\end{bmatrix}\,.
	\end{align}
	\end{subequations}
Applying the conditional Gaussian formula yields
\begin{equation}\label{predu}
\begin{split}
	 & U^* \mid \vb*\varphi, \vb*\theta, \vb*Y = \vb y, \vb*R = \vb f \\
\sim~~ &  
	 \mathcal{N}\left(  \widehat{\vb* k}^{\vb* \varphi}(s^*)\trp (\widehat{\vb K}^{\vb* \varphi})^{-1}\begin{pmatrix}
		\vb y \\ g^{\vb*\theta}(\tilde{\vb u} )+ \vb f 
	\end{pmatrix},
	\kappa(s^*, s^*) - \widehat{\vb* k}^{\vb* \varphi}(s^*)\trp (\widehat{\vb K}^{\vb* \varphi})^{-1} \widehat{\vb* k}^{\vb* \varphi}(s^*)
	\right)\,.
\end{split}
\end{equation}
Note that the Gaussian process predictive mean is a linear combination of kernel functions and their operator-transformed counterparts, and therefore inherits their regularity. Thus, the choice of covariance kernel also determines the smoothness of the resulting solution approximation.

\paragraph{Special case: PDE-constrained Gaussian process} (see, e.g., \cite{raissi2017GPsforlinearPDEs}) Consider a linear partial differential equation
\begin{equation}
	\mathscr{L} u (x)  = f(x)\,,\quad x \in D \subset \mathbb{R}^d
\end{equation}
without any parametrization. In this case, $\mathscr{L}$ is a linear differential operator over a spatial domain $D \subset \mathbb{R}^d$; the coordinate $s = x$ represents the spatial location; and $\vb f$ collects the values of the given forcing at $m'$ locations in $\mathcal{S}'$. The predictive distribution for the test variable $U^* = u(x^*)$ \eqref{predu} then simplifies to
\begin{equation} \label{PDEposterior}
	U^* \mid  \vb*Y = \vb y , \vb*R = \vb f ~\sim \mathcal{N}\left(   \widehat{\vb* k}(x^*)\trp \widehat{\vb K}^{-1}\begin{pmatrix}
		\vb y \\ \vb f 
	\end{pmatrix},
	\kappa(x^*, x^*) - \widehat{\vb* k}(x^*)\trp \widehat{\vb K}^{-1} \widehat{\vb* k}(x^*)
	\right)\,.
\end{equation}
in which $\widehat{\vb* k}$ and $\widehat{\vb K}$ no longer depend on $\vb*\varphi$, as the differential operator is not parametrized. This result can be interpreted as nonparametric Bayesian inference based on both state and derivative observations.

\medskip
These two scenarios demonstrate how the unified Bayesian formulation can be adapted to different inferential goals, namely parameter estimation and equation solving, while relying on the same underlying probabilistic structure and differing only in how inference and marginalization are carried out.

\section{Generality}

The unified Bayesian formulation developed in this work is sufficiently general to encompass a wide range of existing approaches to Gaussian process approximation of differential equations. In particular, many methods that exploit the preservation of Gaussianity under linear differential operators can be interpreted as special cases of the framework introduced in Section 2. In this section, we illustrate this generality by revisiting several representative classes of methods and showing how they can be understood within the same probabilistic structure.

\paragraph{Physics-informed Gaussian processes}

A broad class of approaches, often referred to as physics-informed Gaussian processes \cite{Simo2011,raissi2017GPsforlinearPDEs,raissi2018gppdes,pfortner2022physics,wang2022physics,bai2024gaussian,Ye2025,daniels2025uncertainty}, incorporates differential equation constraints into Gaussian process models by applying linear differential operators directly to the kernels. Such kernel transformations rely on the fact that Gaussian processes are closed under linear operations. From the perspective of the unified Bayesian formulation developed in Section 2, these methods can be understood as instances of the likelihood construction in which observations of the differential operator applied to the solution are incorporated through the variable $\vb*R$. In particular, the use of operator-transformed kernels corresponds to modeling the joint distribution of the solution and its derivatives, while the enforcement of the differential equation is achieved through a matching condition between the `internal force' $\vb*R$ and the external forcing term.

Many existing works follow this general paradigm, including approaches for forward and inverse problems. Depending on the inferential objective, these methods can be interpreted within the two scenarios introduced in Section 3. For forward problems \cite{Simo2011,raissi2017GPsforlinearPDEs,raissi2018gppdes,pfortner2022physics,daniels2025uncertainty}, they typically reduce to Gaussian process regression with both state and derivative observations, corresponding to Scenario II. For inverse problems and parameter estimation \cite{raissi2017GPsforlinearPDEs,bai2024gaussian,Ye2025}, they align with Scenario I, where unknown parameters are inferred while treating the solution as latent.

In particular, \cite{raissi2017GPsforlinearPDEs} formulates parameter estimation as a maximum marginal likelihood problem based on Gaussian process models with operator-transformed kernels, which gives a point estimate through empirical Bayes. The unified Bayesian framework presented in this work instead allows for a fully probabilistic treatment and enables consistent uncertainty quantification for both parameters and solutions. Moreover, recent developments combine Gaussian processes with deep learning to improve scalability and expressiveness of nonlinearity, while retaining the operator-transformed kernel structure \cite{wang2022physics,yan2025pde,yan2025physics}.

Another example is Bayesian Augmented Space Learning \cite{hong2023combining}, which reformulates PDEs with uncertain parameters in an augmented space and employs Gaussian process regression to combine observational data with differential equation constraints for probabilistic response analysis. Its operator-based Gaussian process construction shares the mathematical principles underlying Scenario II.

\paragraph{Latent force models}

Latent force models \cite{Alvarez2009,Alvarez2013} adopt an alternative viewpoint in which the solution is expressed through the action of an inverse differential operator (often represented via a Green’s function):
\begin{equation}
	\mathscr{L}^{-1} f = u\,.
\end{equation}
In this setting, a Gaussian process prior is placed on the latent forcing term $f$, and the solution $u$ inherits a Gaussian process structure through the linear inverse operator $\mathscr{L}^{-1}$. This corresponds to placing the prior on a linearly transformed function rather than directly on the solution. In contrast, the formulation in Section 2 places a Gaussian process prior on the solution $u$ and incorporates the differential equation constraint through a likelihood term involving the transformation $\mathscr{L}u$. Despite this difference, both formulations rely on the same principle that Gaussianity is preserved under linearity, and both yield Gaussian process models that couple the solution with the differential operator. Particularly, latent force models are most naturally aligned with Scenario II in Section 3, where the Gaussian process is used directly for approximating the solution field and the effect of the differential operator is reflected in the covariance structure.

From a nonparametric Bayesian perspective, the two formulations can be viewed as dual constructions: in Scenario II, a Gaussian process prior is placed on the solution $u$ and the likelihood is constructed via the differential operator $\mathscr{L}$, whereas latent force models place the prior on the forcing term $f$ and recover $u$ through the inverse operator $\mathscr{L}^{-1}$. In this sense, they can be regarded as symmetric realizations of the same underlying principle.

\paragraph{Parameter estimation and model identification for dynamical systems}

A large class of problems concerns the estimation of unknown parameters and the identification of governing equations in dynamical systems. Consider, in particular, a system of ordinary differential equations of the form
\begin{equation}
	\frac{\text{d}}{\text{d}t} \vb* q(t) = \vb* g^{\vb*\theta}(\vb*q(t))\,,
\end{equation}
where the parameter vector $\vb*\theta$ characterizes the right-hand side of the system. This setting can be viewed as a multidimensional generalization of the special case discussed in Scenario I in Section 3. This extension to coupled nonlinear ODE systems is further illustrated by the Lotka--Volterra example in Appendix A.

From the perspective of the unified Bayesian formulation, these problems are treated by placing a Gaussian process prior on the solution trajectory and inferring the parameters $\vb*\theta$ through marginalization over the latent state. The resulting posterior distribution provides a probabilistic characterization of the parameters, conditioned on observed data and the dynamical structure.

Different methodologies arise depending on how the parameter vector $\vb*\theta$ is interpreted. In classical parameter estimation, $\vb*\theta$ represents unknown physical parameters governing the system, and the task reduces to Bayesian inference in dynamical systems based on a Gaussian process prior over the trajectory \cite{yang2021inference,ye2024gaussian}. In model identification problems, $\vb*\theta$ parametrizes the form of the right-hand side. For example, in sparse identification of nonlinear dynamics (SINDy) \cite{brunton2016discovering}, $\vb*\theta$ corresponds to coefficients in a prescribed dictionary of candidate functions; in neural ordinary differential equations \cite{Chen2018}, $\vb*\theta$ represents the parameters of a neural network defining $g^{\vb*\theta}$; and in data-driven model reduction \cite{ghattas2021acta,Benjamin2016,McQuarrie2025}, $\vb*\theta$ encodes the structure of a reduced-order dynamical system.

Despite these different interpretations, all these techniques can be understood within the same inferential framework, in which the parameters $\vb*\theta$ are inferred while the solution is treated as latent. In this sense, parameter estimation and model identification for dynamical systems constitute variations of the same underlying problem within the unified Bayesian formulation.

\paragraph{Extension to weak form}

The unified Bayesian framework presented above also extends to weak formulations of differential problems, subject to modifications analogous to the treatment for the strong form \eqref{DE}. The weak form of \eqref{DE} can typically be written as follows: find $u\in \mathcal{U}$, where $\mathcal{U}$ is a suitable Banach space of admissible solution functions, such that
\begin{equation}\label{weakform}
	a^{\vb*\varphi}(u, v) - \langle \mathfrak{g}^{\vb*\theta}(u), v \rangle_{\mathcal{V}'\times \mathcal{V}}= l(v)\,,\quad \forall v \in \mathcal{V}\,, 
\end{equation}
where $\mathcal{V}$ is the test (Banach) space, $a^{\vb*\varphi}(\cdot, \cdot): \mathcal{U}\times \mathcal{V} \to \mathbb{R}$ is a bilinear form induced by the differential operator $\mathscr{L}^{\vb*\varphi}$, $\langle\cdot,  \cdot \rangle_{\mathcal{V}'\times \mathcal{V}}$ denotes the duality pairing between $\mathcal{V}'$ and $\mathcal{V}$, $\mathfrak{g}^{\vb*\theta} : \mathcal{U} \to \mathcal{V}'$ is a nonlinear operator corresponding to the term $g^{\vb*\theta}$ in \eqref{DE}, and $l\in \mathcal{V}'$ is a linear functional corresponding to the forcing term $f$. This setting can be incorporated into the unified Bayesian framework, with several modifications to the interpretation of the observables. 

While the direct observations $\vb* Y = \vb y$ at $m$ locations in $\mathcal{S}\subset D$ remain unchanged, the interpretation of the observable $\vb*R$ and the `derivatives' $\vb*Z$ must be adapted in an analogous manner. Specifically, $\vb* R$ now represents the functional $a^{\vb*\varphi}(u, \cdot) - \langle \mathfrak{g}^{\vb*\theta}(u), \cdot \rangle_{\mathcal{V}'\times \mathcal{V}}$ evaluated at a set of $m'$ independent test functions $\vb*\phi = ( \phi_1, \cdots, \phi_{m'}) \trp $ all belonging to $\mathcal{V}$, while $\vb* Z$ corresponds to the functional $a^{\vb*\varphi}(u, \cdot)$ evaluated at the same set of test functions. Analogous to the strong-form construction, $\vb* R$ is approximated by $\vb*Z - \langle \mathfrak{g}^{\vb*\theta}(\tilde{u}), \Phi \rangle_{\mathcal{V}'\times \mathcal{V}}$, where $\tilde{u}$ denotes a data-driven estimate of $u$, for example $\tilde{u}(\cdot) = \kappa(\cdot, \mathcal{S})(\vb K + \sigma_y^2\vb I)^{-1}\vb y$, and the approximation error is modeled as white noise with variance $\sigma_f^2$. Moreover, the observed data for $\vb*R$ is now given by $\vb f = (l(\phi_1),\cdots, l(\phi_{m'}))\trp$. 

With these modifications, the unified Bayesian formulation and the two typical scenarios remain applicable. The main difference lies in the operator-transformed kernel terms, which now take the following form:
\begin{subequations} \label{wfterms}
	\begin{align}
		(\mathscr{L}^{\vb* \varphi}\otimes \text{Id})\kappa(\mathcal{S}',\mathcal{S}) ~& \leftarrow (a^{\vb* \varphi}(\cdot, \vb*\phi)\otimes \text{Id})\kappa(\cdot,\mathcal{S}) \,, \\
		(\text{Id} \otimes \mathscr{L}^{\vb* \varphi} )\kappa(\mathcal{S},\mathcal{S}') ~& \leftarrow (\text{Id} \otimes a^{\vb* \varphi} (\cdot, \vb*\phi))\kappa(\mathcal{S},\cdot) \,, \\
		(\mathscr{L}^{\vb* \varphi}\otimes \mathscr{L}^{\vb* \varphi})\kappa(\mathcal{S}',\mathcal{S}') ~& \leftarrow (a^{\vb* \varphi}(\cdot, \vb*\phi) \otimes a^{\vb* \varphi}(\cdot, \vb*\phi))\kappa(\cdot,\cdot) \,,
	\end{align}
\end{subequations}
that is, the strong-form collocation is replaced by weak-form functional evaluation. To clarify the notation, the matrix components are given by $[(a^{\vb* \varphi}(\cdot, \vb*\phi)\otimes \text{Id})\kappa(\cdot,\mathcal{S})]_{ij} = a^{\vb* \varphi}(\kappa(\cdot, s_j), \phi_i)$, $[(\text{Id} \otimes a^{\vb* \varphi} (\cdot, \vb*\phi))\kappa(\mathcal{S},\cdot)]_{ij} = a^{\vb* \varphi}(\kappa(\cdot, s_i), \phi_j)$, and $[(a^{\vb* \varphi}(\cdot, \vb*\phi) \otimes a^{\vb* \varphi}(\cdot, \vb*\phi))\kappa(\cdot,\cdot)]_{ij} = (a^{\vb* \varphi}(\cdot, \phi_i) \otimes a^{\vb* \varphi}(\cdot, \phi_j))\kappa(\cdot,\cdot)$. Note that we here assume $\kappa(\cdot, s) = \kappa(s, \cdot) \in \mathcal{U}$ for all $s\in D$, so that the kernel functions are admissible in the weak form. This shows that the unified Bayesian framework extends naturally to variational settings that are commonly used in numerical analysis.

Appendix B further illustrates the weak-form formulation through an analytical derivation of the Gaussian process posterior, and also highlights its connections to Bayesian linear regression and the classical Petrov-Galerkin method.

\section{Concluding remarks}

This note presents a unified Bayesian perspective on Gaussian process approximation of differential equations. Two typical scenarios, namely parameter estimation and solution approximation, are shown to arise naturally, and a range of existing relevant methods can be interpreted within this framework. More broadly speaking, this unified view contributes to a necessary consolidation of the rapidly growing literature in this area. By clarifying common structures, connections, and relationships, it provides a foundation for more systematic developments in Gaussian process-based data-driven methods for differential equations.

Looking ahead, a primary challenge is the curse of dimensionality, which remains a fundamental limitation when applying Gaussian processes to high-dimensional partial differential equations. Addressing this will likely require an \textit{effective} integration of Gaussian processes and deep learning, along with the development of corresponding theoretical foundations in statistical learning and numerical analysis.

\section*{Acknowledgment}

This work was partially supported by the Wallenberg AI, Autonomous Systems and Software Program (WASP) funded by the Knut and Alice Wallenberg Foundation. The author also acknowledges financial support from the Swedish Research Council (VR) under No. 2025-04911 and from The Crafoord Foundation under No. 20250657. The author expresses his appreciation to Dr. Anirban Chaudhuri (UT Austin) and Dr. Shane McQuarrie (BYU) for fruitful relevant discussions, and to Dr. Sarah Vollert (Umeå) for her valuable feedback.

\section*{Declaration of AI use}

During the preparation of this work the author used GPT-5 and GPT-6 for language polishing, particularly to improve textual clarity, and for assistance with LaTeX formatting. After using these tools, the author reviewed and edited the content as needed and takes full responsibility for the content of the published article. All technical and mathematical contributions are the author's own.

\section*{Appendix A. Parameter estimation for the Lotka--Volterra system}

Consider the Lotka--Volterra system
\begin{equation*}
\begin{cases}
\frac{\text{d} x_1}{\text{d}t} (t) =\alpha x_1(t)-\beta x_1(t)x_2(t),\\
\frac{\text{d} x_2}{\text{d}t} (t) =\delta x_1(t)x_2(t)-\gamma x_2(t),
\end{cases}
\end{equation*}
where $\vb*\theta=(\alpha,\beta,\gamma,\delta)\trp$ denotes the unknown governing parameters to infer. Let $\mathcal{S} =\{t_1,\ldots,t_m\}$ denote the observation times, and $\mathcal{S}'=\{t'_1,\cdots,t'_{m'}\}$ the time instances used for derivative matching. Given noisy observations $\vb y_i$ of each state $x_i$ at $\mathcal{S}$, independent Gaussian process priors $x_i\sim\mathcal{GP}(0,\kappa_i)$ are assigned, $i=1,2$.

Define
\begin{align*}
\vb K_i & =\kappa_i(\mathcal S, \mathcal S),  & 
\vb C_i & =\kappa_i(\mathcal S', \mathcal S), \\
\vb D_i & = \left (\tfrac{\text{d}}{\text{d}t}\otimes\text{Id} \right)\kappa_i(\mathcal S', \mathcal S), &  \text{and} \quad
\vb H_i & =\left(\tfrac{\text{d}}{\text{d}t} \otimes \tfrac{\text{d}}{\text{d}t} \right)\kappa_i(\mathcal S',\mathcal S')\,, \quad i = 1,2.
\end{align*}
The Gaussian process regression for each state component gives the predictive means of the states and their derivatives at $\mathcal S'$, together with the covariance matrices for the time-derivative estimates:
\begin{equation*}
\begin{split}
\tilde{\vb x}_i =\vb C_i (\vb K_i +\sigma_{y,i}^2\vb I_m) ^{-1}\vb y_i,\quad 
\widehat{\vb* \mu}_{z,i} & =\vb D_i (\vb K_i +\sigma_{y,i}^2\vb I_m) ^{-1} \vb y_i, \\
\text{and} \quad 
\widehat{\vb* \Sigma}_{z,i} & =\vb H_i -\vb D_i (\vb K_i + \sigma_{y,i}^2\vb I_m) ^{-1} \vb D_i \trp +\sigma_{f,i}^2\vb I_{m'}\,.
\end{split}
\end{equation*}
Although the ODE system is nonlinear in its states, it is linear in the governing parameters $\vb*\theta$. Using the Gaussian process predictive means to approximate the nonlinear terms, define
\begin{equation*}
\vb G_1=
\begin{bmatrix}
\tilde{\vb x}_1&
-\tilde{\vb x}_1\odot\tilde{\vb x}_2
\end{bmatrix} \quad \text{and} \quad 
\vb G_2=
\begin{bmatrix}
-\widetilde{\vb x}_2&
\tilde{\vb x}_1\odot\tilde{\vb x}_2
\end{bmatrix},
\end{equation*}
in which $\odot$ denotes elementwise multiplication, and
\begin{equation*}
\widehat{\vb* \mu}_{z}
=\begin{bmatrix} \widehat{\vb* \mu}_{z,1}\\ \widehat{\vb* \mu}_{z,2} \end{bmatrix},
\quad
\vb G=\text{blkdiag}(\vb G_1,\vb G_2), \quad \text{and} \quad
\widehat{\vb* \Sigma}_{z} = \text{blkdiag}(\widehat{\vb* \Sigma}_{z,1},\widehat{\vb* \Sigma}_{z,2}).
\end{equation*}

Hence, the likelihood function \eqref{like_f} given by derivative matching is
\begin{equation*}
\pi_\text{like}(\vb y_1, \vb y_2, \vb f = \vb 0 \mid \vb*\theta) 
\propto
\exp\left[ -\frac12 (\vb G\vb*\theta- \widehat{\vb* \mu}_{z})\trp
{\widehat{\vb* \Sigma}_{z}}^{-1} (\vb G\vb*\theta- \widehat{\vb* \mu}_{z}) \right].
\end{equation*}
With a non-degenerate Gaussian prior
\begin{equation*}
\pi_\text{prior} (\vb* \theta) \propto
\exp\left[ -\frac12 (\vb*\theta- \vb*\theta_0)\trp
{\vb* \Lambda}^{-1} (\vb*\theta- \vb*\theta_0) \right],
\end{equation*}
the resulting posterior takes the closed form
\begin{equation*}
\begin{split}
& \pi_\text{post} (\vb* \theta \mid \vb y_1, \vb y_2, \vb f = \vb 0 )
 = \mathcal N(\vb*\theta \mid \vb*\mu_\theta,\vb*\Sigma_\theta)\,,\\
 \text{with} \quad 
& \vb*\Sigma_\theta =\left( \vb G \trp {\widehat{\vb* \Sigma}_{z}} ^{-1}\vb G +\vb*\Lambda^{-1} \right)^{-1} \quad \text{and} \quad
\vb*\mu_\theta = \vb*\theta_0  + \vb*\Sigma_\theta
 \vb G \trp {\widehat{\vb* \Sigma}_{z}} ^{-1} \left( \widehat{\vb* \mu}_{z} - \vb G \vb*\theta_0 \right).
 \end{split}
\end{equation*}
This example illustrates how the derivative-matching formulation extends to coupled nonlinear ODE systems, in which Gaussian process regression is performed using observations at $\mathcal{S}$, while parameter estimation is based on the approximate differential equation constraints evaluated at $\mathcal{S}'$. The linear dependence on the system parameters yields a closed-form Gaussian posterior. The result here is consistent with that in \cite{ye2024gaussian}.

\section*{Appendix B. Solving a PDE in its weak form}

Consider a special case of the weak formulation \eqref{weakform}: find $u\in\mathcal U$ such that
\begin{equation*}
    a(u,v)= l(v),\qquad \forall v\in\mathcal V.
\end{equation*}
In this case, the aim is to solve the equation without observational data $\vb{y}$.
As discussed in Section 4 ``Extension to weak form'', take $m'$ linearly independent test functions
$\vb*{\phi}=(\phi_1,\ldots,\phi_{m'})\trp\in\mathcal V^{m'}$
and define the forcing data
\begin{equation*}
    \vb{f}=\bigl(l(\phi_1),\ldots,l(\phi_{m'})\bigr)\trp.
\end{equation*}
Using \eqref{wfterms} and the PDE-constrained Gaussian process posterior of the solution \eqref{PDEposterior}, the conditional distribution for $u(x^*)$ is
\begin{equation*}
    U^*\mid \vb*R=\vb{f} ~\sim\mathcal N\!\left(
        \widehat{\vb*{k}}(x^*)\trp\widehat{\vb{K}}^{-1}\vb{f},\;
        \kappa(x^*,x^*)-\widehat{\vb*{k}}(x^*)\trp
        \widehat{\vb{K}}^{-1}\widehat{\vb*{k}}(x^*)
    \right),
\end{equation*}
in which
\begin{equation*}
    \widehat{\vb*{k}}(x^*)
        =\bigl(a(\cdot,\vb*{\phi})\otimes\mathrm{Id}\bigr)
          \kappa(\cdot,x^*) \qquad \text{and} \qquad 
    \widehat{\vb{K}}
        =\bigl(a(\cdot,\vb*{\phi})\otimes
          a(\cdot,\vb*{\phi})\bigr)\kappa(\cdot,\cdot)
          +\sigma_f^2\vb{I}_{m'}.
\end{equation*}

Next, consider a finite-dimensional approximation of the solution in terms of
$N$ basis functions $\vb*{\psi}=(\psi_1,\ldots,\psi_N)\trp\in\mathcal U^N$:
\begin{equation*}
    u(x)\approx\vb*{w}\trp\vb*{\psi}(x).
\end{equation*}
Assigning a Gaussian prior $\vb*{w}\sim\mathcal N(\vb*{0},\vb*{\Gamma})$ to the combination coefficients $\vb* w$,
    with $\vb*{\Gamma}\succ \vb 0$, induces a zero-mean Gaussian process prior on the approximate
solution with finite-rank covariance kernel
\begin{equation*}
    \kappa(x,x')=\vb*{\psi}(x)\trp
      \vb*{\Gamma}\vb*{\psi}(x').
\end{equation*}
Substituting this kernel into the expressions above gives
\begin{equation*}
    \widehat{\vb*{k}}(x^*)=\vb{A}\vb*{\Gamma}\vb*{\psi}(x^*) \qquad \text{and} \qquad
    \widehat{\vb{K}}=\vb{A}\vb*{\Gamma}\vb{A}\trp
        +\sigma_f^2\vb{I}_{m'},
\end{equation*}
where $\vb{A}\in\mathbb R^{m'\times N}$ has entries
$\vb A_{ij}=a(\psi_j,\phi_i)$.
Consequently, the predictive mean and variance of solution $U^* = u(x^*)$ can be expressed respectively as
\begin{equation*}
    \mathbb E[U^*\mid \vb*R=\vb{f}] =
        \vb*{\psi}(x^*)\trp\vb*{\mu}_w \qquad \text{and} \qquad
    \mathbb Var[U^*\mid \vb*R=\vb{f}] =
        \vb*{\psi}(x^*)\trp\vb*{\Sigma}_w
        \vb*{\psi}(x^*)\,,
\end{equation*}
with $\vb*{\mu}_w =\vb*{\Gamma}\vb{A}\trp \bigl(\vb{A}\vb*{\Gamma}\vb{A}\trp +\sigma_f^2\vb{I}_{m'}\bigr)^{-1}\vb{f}$, and $\vb*{\Sigma}_w =\vb*{\Gamma}-\vb*{\Gamma}\vb{A}\trp\bigl(\vb{A}\vb*{\Gamma}\vb{A}\trp +\sigma_f^2\vb{I}_{m'}\bigr)^{-1}\vb{A}\vb*{\Gamma}$. Equivalently,
\begin{equation*}
    \vb*{\Sigma}_w
        =\left(\sigma_f^{-2}\vb{A}\trp\vb{A}
            +\vb*{\Gamma}^{-1}\right)^{-1}
    \qquad \text{and} \qquad
    \vb*{\mu}_w
        = \sigma_f^{-2}\vb*{\Sigma}_w
          \vb{A}\trp\vb{f}.
\end{equation*}
Thus, the coefficient vector has the posterior distribution
$\vb*{w}\mid \vb* R=\vb{f}\sim\mathcal N(\vb*{\mu}_w,\vb*{\Sigma}_w)$.
Its density can be factored as
\begin{equation*}
    \exp\left[-\frac12
        (\vb*{w}-\vb*{\mu}_w)\trp
        \vb*{\Sigma}_w^{-1}
        (\vb*{w}-\vb*{\mu}_w)\right]
    \propto \exp\left[-\frac{1}{2\sigma_f^2}
        \|\vb{A}\vb*{w}-\vb{f}\|_2^2\right]
        \exp\left[-\frac12\vb*{w}\trp
        \vb*{\Gamma}^{-1}\vb*{w}\right].
\end{equation*}
The two factors correspond to the likelihood
$\vb{f}-\vb{A}\vb*{w}\sim\mathcal N(\vb*{0},\sigma_f^2\vb{I}_{m'})$
and the prior $\vb*{w}\sim\mathcal N(\vb*{0},\vb*{\Gamma})$,
respectively, as in Bayesian linear regression.
The maximum a posteriori (MAP) estimate therefore solves
\begin{equation*}
    \underset{\vb*{w}\in\mathbb R^N}{\min}
      \left\{\frac{1}{\sigma_f^2}
        \|\vb{A}\vb*{w}-\vb{f}\|_2^2
        +\|\vb*{w}\|_{\vb*{\Gamma}^{-1}}^2\right\}
    =\underset{\vb*{w}\in\mathbb R^N}{\min}
      \left\{\frac{1}{\sigma_f^2}\sum_{i=1}^{m'}
        \bigl[a(\vb*{w}\trp\vb*{\psi},\phi_i)
              - l(\phi_i)\bigr]^2
        +\|\vb*{w}\|_{\vb*{\Gamma}^{-1}}^2\right\}.
\end{equation*}
In the limit of vanishing prior precision (uninformative prior),
$\vb*{\Gamma}^{-1}\to\vb*{0}$, the Tikhonov regularization term vanishes; additionally 
if $N=m'$ and $\vb{A}$ is nonsingular, the MAP equation reduces to
the Petrov-Galerkin scheme: find
$u_h\in\text{span}\{\psi_1,\ldots,\psi_N\}$ such that
\begin{equation*}
    a(u_h,\phi_i)=l(\phi_i),\qquad i=1,\cdots,N.
\end{equation*}

\bibliographystyle{abbrv}
\bibliography{refs.bib}

\end{document}